\documentclass{ifacconf}

\usepackage{amsmath} 
\usepackage{amssymb}
\usepackage{amsfonts}
\usepackage{dsfont}
\usepackage{natbib}        % required for bibliography

\newcommand{\bq}{\begin{equation}}
\newcommand{\eq}{\end{equation}}
\newcommand{\mR}{\mathbb{R}}
\newcommand{\mP}{\mathds{P}}

\newcommand{\A}{\mathcal{A}}

\newcommand{\mL}{\mathbb{L}}

\newcommand{\T}{\mathcal{T}}

\newcommand{\cL}{\mathcal{L}}
\newcommand{\Q}{\mathcal{Q}}
\newcommand{\Exp}{\mathrm{Exp}\,}

\newtheorem{df}[thm]{Definition}
\newtheorem{ex}[thm]{Example}
\begin{document}
\begin{frontmatter}

\title{Towards Control by Interconnection of Port-Thermodynamic Systems} 

%\thanks[footnoteinfo]{Sponsor and financial support acknowledgment
%goes here. Paper titles should be written in uppercase and lowercase
%letters, not all uppercase.}
%\thanks[footnoteinfo]{This paper was written as the result of the stay in 2017 of the first author as invited professor with the research group DYCOP of the Laboratoire d'Automatique et de G\'enie des Proc\'ed\'es (LAGEP), Universit\'e Claude Bernard Lyon 1.}

\author{Arjan van der Schaft} 
%\author[Second]{Bernhard Maschke} 

\address{Bernoulli Institute for Mathematics, Computer Science and AI, \\
Jan C. Willems Center for Systems and Control, \\University of Groningen, the Netherlands \\(e-mail: a.j.van.der.schaft@rug.nl)}
%\address[Second]{Universit\'e Claude Bernard Lyon 1, CNRS, LAGEP, France \\(e-mail: bernhard.maschke@univ-lyon1.fr)}

\begin{abstract}                % Abstract of not more than 250 words.
The power conserving interconnection of port-thermodynamic systems via their power ports results in another port-thermodynamic system, while the same holds for the rate of entropy increasing interconnection via their entropy flow ports. Control by interconnection of port-thermodynamic systems seeks to control a plant port-thermodynamic system by the interconnection with a controller port-thermodynamic system. The stability of the interconnected port-thermodynamic system is investigated by Lyapunov functions based on generating functions for the submanifold characterizing the state  properties as well as additional conserved quantities. Crucial tool is the use of canonical point transformations on the symplectized thermodynamic phase space.
\end{abstract}


\end{frontmatter}
%===============================================================================

\section{INTRODUCTION}
Since the 1970s, based on Gibbs' fundamental thermodynamic relation, {\it contact geometry} has been recognized as an appropriate geometric framework for macroscopic thermodynamics; see e.g. \cite{hermann, mrugala1, mrugala2, mrugala4, eberard,favache09, favache10, bravetti17, bravetti, deleon, hudon, ramirez3, gromov}. In \cite{balian} it was argued that the distinction between the contact-geometric description of the {\it energy} and {\it entropy} representation of thermodynamic systems can be resolved by {\it symplectization} of contact manifolds; a concept which is known in differential geometry, cf. \cite{arnold, libermann}. Subsequently in \cite{valparaiso1, entropy, valparaiso2, jgp} this viewpoint was extended, leading to the general definition of {\it port-thermodynamic systems}. 

In the present paper we will initiate a methodology of {\it control by interconnection} for port-thermodynamic systems. In this approach we seek to control a given (plant) port-thermodynamic system by interconnecting it with another (controller) port-thermodynamic system. Following \cite{entropy} the interconnection of port-thermodynamic systems via {\it power ports} or {\it entropy flow ports} results in an interconnected system that is again a port-thermodynamic system. Our main results concern {\it regulation} of port-thermodynamic systems, where the plant port-thermodynamic system is sought to be asymptotically stabilized at a desired set-point. This is achieved by the construction of a Lyapunov function, employing, next to the generating function of the Liouville submanifold describing the state properties, the presence of conserved quantities. The approach mimics control by interconnection of port-Hamiltonian systems; however with some fundamental differences as discussed in the Conclusions Section \ref{sec:concl}. 

The paper provides in Section \ref{sec:recall} a recall of the definition of port-thermodynamic systems from \cite{entropy, jgp}. Section \ref{sec:asymp} deals with (asymptotic) stabilization of port-thermodynamic systems using conserved quantities and interconnection with damper systems, while Section \ref{sec:int} initiates the general control by interconnection methodology.

\section{RECALL OF PORT-THERMODYNAMIC SYSTEMS}\label{sec:recall}
Consider a simple thermodynamic system, such as a single gas in a compartment with volume $V$ and pressure $P$ at temperature $T$. It is well-known that the {\it state properties} of the gas are described by a $2$-dimensional submanifold of the ambient space $\mR^5$ (the {\it thermodynamic phase space}) with coordinates $E$ (energy), $S$ (entropy), $V$, $P$, and $T$. Such a submanifold characterizes the properties of the gas (e.g., an ideal gas, or a Van der Waals gas), and all of them share the following property. Define the {\it Gibbs one-form} on the thermodynamic phase space $\mR^5$ as
\bq
\label{Gibbs}
\theta:=dE - TdS +PdV
\eq
Then $\theta$ is {\it zero} restricted to the submanifold characterizing the state properties. This is called {\it Gibbs' fundamental thermodynamic relation}. Geometrically the Gibbs one-form $\theta$ defines a {\it contact form} on $\mR^5$, and any submanifold $L$ capturing the state properties of the thermodynamic system is a submanifold of maximal dimension restricted to which the contact form $\theta$ is zero. Such submanifolds are called {\it Legendre submanifolds} of the {\it contact manifold} $(\mR^5, \theta)$.

By expressing the extensive variable $E$ as a function $E=\bar{E}(S,V)$ of the two remaining extensive variables $S$ and $V$, Gibbs' fundamental relation implies that the Legendre submanifold $L$ specifying the state properties is given as
\bq
\label{L1}
L=\{(E,S,V,T,P) \mid E=\bar{E}(S,V), T= \frac{\partial \bar{E}}{\partial S}, -P= \frac{\partial \bar{E}}{\partial V} \}
\eq
Hence $L$ is completely described by the energy function $\bar{E}(S,V)$, whence the name {\it energy representation} for \eqref{L1}. 
{\it Another} option for describing $L$ is the {\it entropy representation}. This option is motivated from a modeling point of view by noting that often thermodynamic systems are formulated by first listing the {\it balance laws} for all the extensive variables except for the entropy $S$, and then expressing $S$ as a function $S=\bar{S}(E,V)$. This leads to the representation of $L \subset \mR^5$, given as
\bq
\label{L3}
L:= \{(E,S,V,T,P) \mid S=\bar{S}(E,V), \frac{1}{T}= \frac{\partial \bar{S}}{\partial E}, \frac{P}{T}= \frac{\partial \bar{S}}{\partial V} \}
\eq
Geometrically the entropy representation corresponds to the {\it modified} Gibbs contact form
\bq
\label{Gibbs1}
\widetilde{\theta}:=dS -\frac{1}{T}dE - \frac{P}{T}dV,
\eq
which is obtained from the original Gibbs contact form $\theta$ in \eqref{Gibbs} by division by $-T$ (called {\it conformal equivalence}). In this way the Gibbs fundamental relation is rewritten as $\widetilde{\theta}|_L=0$, and the intensive variables become $\frac{1}{T},\frac{P}{T}$.

As argued in \cite{entropy,jgp}, continuing on \cite{balian}, the contact-geometric view on thermodynamics has two shortcomings:\\
(1) Switching from the energy representation $E=\bar{E}(S,V)$ to the entropy representation $S=\bar{S}(E,V)$ corresponds to replacing the Gibbs form $\theta$ by the modified Gibbs form $\widetilde{\theta}$, and thus leads to a {\it different} (although conformally equivalent) contact-geometric description.\\
(2) The contact-geometric description does not make a clear distinction between, on the one hand, the extensive variables $E,S,V$ and, on the other hand, the intensive variables $T,-P$ (energy representation), or $\frac{1}{T},\frac{P}{T}$ (entropy representation). 

The way to remedy these shortcomings is to {\it extend} the contact manifold by one extra dimension to a symplectic manifold, in fact a {\it cotangent bundle}, with an additional {\it homogeneity} structure. This construction is rather well-known in differential geometry \cite{arnold, libermann}, but was advocated within a thermodynamics context only in \cite{balian}, and then followed up in \cite{valparaiso1, entropy, valparaiso2, jgp}. As argued in \cite{entropy, jgp} this point of view has computational advantages as well.

For a simple thermodynamic system with extensive variables $E,S,V$ and intensive variables $T,-P$, the construction amounts to replacing the intensive variables $T,-P$ by their {\it homogeneous coordinates} $p_E,p_S,p_V$ with $p_E \neq 0$, i.e., 
\bq
T= \frac{p_S}{-p_E}, \; -P= \frac{p_V}{-p_E}
\eq
Equivalently, the intensive variables $\frac{1}{T}, \frac{P}{T}$ in the {\it entropy representation} are represented as
\bq
\frac{1}{T} = \frac{p_E}{-p_S}, \; \frac{P}{T}= \frac{p_V}{-p_S},
\eq
where $p_S \neq 0$. This means that the {\it two} contact forms $\theta=dE - TdS +PdV$ and $\widetilde{\theta}=dS -\frac{1}{T}dE - \frac{P}{T}dV$ are replaced by a {\it single} symmetric expression
\bq
\alpha:=p_EdE + p_SdS + p_VdV,
\eq
The one-form $\alpha$ is the canonical {\it Liouville one-form} on the cotangent bundle $T^*\mR^3$, with $\mR^3$ the space of extensive variables $E,S,V$. Thus the thermodynamic phase space $\mR^5$ has been replaced by $T^*\mR^3$. More precisely, by definition of homogeneous coordinates the vector $(p_E,p_S,p_V)$ is different from the zero vector, and hence the space with coordinates $E,S,V,p_E,p_S,p_V$ is actually the cotangent bundle $T^*\mR^3$ {\it minus} its zero section; denoted as $\T^*\mR^3$.

Any $2$-dimensional Legendre submanifold $L \subset \mR^5$ describing the state properties is now replaced by a $3$-dimensional submanifold $\cL \subset \T^*\mR^3$, given as
\bq
\cL=\{(E,S,V,p_E,p_S,p_V) \mid (E,S,V, \frac{p_S}{-p_E}, \frac{p_V}{-p_E}) \in L  \}
\eq
It turns out that $\cL$ is a {\it Lagrangian submanifold} of $\T^*\mR^3$ with symplectic form $\omega:=d\alpha$, with an additional property of {\it homogeneity}. Namely, whenever $(E,S,V,p_E,p_S,p_V) \in \cL$, then also $(E,S,V,\lambda p_E,\lambda p_S, \lambda p_V) \in \cL$, for any non-zero $\lambda \in \mR$. Such Lagrangian submanifolds are shown \cite{entropy,jgp} to be fully characterized as maximal manifolds restricted to which the Liouville one-form $\alpha=p_EdE + p_SdS + p_VdV$ is zero. They have been called {\it Liouville submanifolds} of $\T^*\mR^3$ in \cite{jgp}.

This is immediately extended to {\it general} thermodynamic phase spaces. For instance, in the case of multiple chemical species the Gibbs form $\theta$ extends to $dE - TdS +PdV - \sum_k \mu_kdN_k$, where $N_k$ and $\mu_k$, $k=1, \cdots,s,$ are the mole numbers, respectively, chemical potentials of the $k$-th species. Correspondingly, the contact manifold $\mR^5 \times \mR^{2s}$ is replaced by the cotangent bundle without zero-section $\T^*\mR^{3 + s}$, with extensive variables $E,S,V,N_1, \cdots, N_s$ and Liouville form
\bq
p_EdE + p_SdS + p_VdV + p_1dN_1 + \cdots + p_sdN_s,
\eq
where $\mu_1 =\frac{p_1}{-p_E}, \cdots, \mu_s =\frac{p_s}{-p_E}$. 

In general, thermodynamic systems will be described on cotangent bundles without zero-section $\T^*\Q$. Here $\Q$ is the $(n+2)$-dimensional manifold of {\it all} extensive variables, denoted by $q^e \in \Q$. We will single out the special extensive variables $E$ (energy) and $S$ (entropy), and write $q^e=(E,S,q)$ with $q$ denoting the {\it remaining} $n$ extensive variables (such as volume and mole numbers). Correspondingly we will denote the coordinates for the cotangent spaces $T_{q^e}^* \Q$ by $p^e=(p_E,p_S,p)$ (called the co-extensive variables).

Next to the {\it state properties} described by a Liouville submanifold $\cL \subset \T^*\Q$ the {\it dynamics} of a thermodynamic system is described by a Hamiltonian vector field $X_K$ on $\T^*\Q$, with the extra requirement that the Hamiltonian $K$ is homogeneous of degree $1$ in the co-extensive variables $p^e$. Equivalently \cite{entropy, jgp} this means that we consider Hamiltonian vector fields that leave the Liouville form $\alpha$ on $\T^*\Q$ invariant. For simplicity of terminology the Hamiltonians $K: \T^*\Q \to \mR$ that are homogeneous of degree $1$ in $p^e$, and their corresponding Hamiltonian vector fields $X_K$, will be simply called {\it homogeneous} in the sequel. 

Furthermore, the homogeneous Hamiltonian vector field $X_K$ should leave the state properties, i.e., the Liouville submanifold $\cL$, invariant. This is equivalent \cite{entropy, jgp, libermann} to the homogeneous Hamiltonian $K$ being zero on $\cL$.
In order to describe the interaction of the thermodynamic system with its surrounding we will split the homogeneous Hamiltonian $K$ into two parts, i.e.,
\bq
K^a + K^s, 
\eq
where $K^a: \T^*\Q^* \to \mR$ is the homogeneous Hamiltonian corresponding to the {\it autonomous} dynamics due to internal non-equilibrium conditions, while $K^s: \T^*\Q^* \times \mR^m \to \mR$, with $u \in \mR^m$ a vector of {\it input} variables, represents the interaction of the system with its surrounding. We assume that $K^s$ for $u=0$ is equal to zero, and thus we may write
\bq
K^s = \sum_{j=1}^m K_j^c u_j,
\eq
for certain functions $K^c_j: \T^*\Q^* \times \mR^m \to \mR$ (interaction or control Hamiltonians). In most situations the dependence of $K^s$ on $u$ is {\it linear}, in which case $K^c_j: \T^*\Q^* \to \mR$; i.e., independent of $u$.
%Recall that $K^a, K^s$ and $K^c_j$ are all homogeneous of degree $1$ in $p$ and zero on $\cL$ for all $u \in \mR^m$. 

By invoking Euler's theorem on homogeneous functions, homogeneity of degree $1$ in $p^e$ means
\bq
\label{KaKc}
\begin{array}{rcll}
K^a &= & p_E\frac{\partial K^a}{\partial p_E} +  p_S\frac{\partial K^a}{\partial p_S} + p_1\frac{\partial K^a}{\partial p_1}+ \cdots + p_n\frac{\partial K^a}{\partial p_n}\\[3mm]
K^c_j & = & p_E\frac{\partial K^c_j}{\partial p_E} +  p_S\frac{\partial K^c_j}{\partial p_S} + p_1\frac{\partial K^c_j}{\partial p_1} +  \cdots + p_n\frac{\partial K^c_j}{\partial p_n} ,
\end{array}
\eq
where all partial derivatives of $K^a$ and $K^c_j$ with respect to $p_E,p_S,p$ are homogeneous of degree $0$ in $p^e=(p_E,p_S,p)$. 

Finally, the class of allowable autonomous Hamiltonians $K^a$ is further restricted by the {\it First and Second Law} of thermodynamics. In fact, since the evolution of $E$ in the autonomous dynamics $X_{K^a}$ arising from non-equilibrium conditions is given by $\dot{E}=\frac{\partial K^a}{\partial p_E}$ the First Law implies that any Hamiltonian $K^a$ should satisfy $\frac{\partial K^a}{\partial p_E}|_{\cL}=0$.
Furthermore, $\dot{S}$ in the autonomous dynamics $X_{K^a}$ is given by $\frac{\partial K^a}{\partial p_S}$. Hence by the Second Law necessarily $\frac{\partial K^a}{\partial p_S}|_{\cL} \geq 0$. 

These two constraints need not hold for the control Hamiltonians $K^c_j$. In fact, the analogous terms in the control Hamiltonians may be utilized to define natural {\it outputs}. First option is to define the output vector as the $m$-dimensional row vector ($p$ for power)
\bq
y_p=\frac{\partial K^c}{\partial p_E},
\eq
where $K^c=(K^c_1,\cdots,K^c_m)$. It follows that along the complete dynamics $X_K$ on $\cL$, with $K=K^a + \sum_{j=1}^m K_j^c u_j$,
\bq
\frac{d}{dt} E= y_pu
\eq
Thus $y_p$ is the vector of {\it power conjugate} outputs corresponding to the input vector $u$. We call the pair $(u, y_p)$ the {\it power port} of the system. 

Similarly, by defining the output vector as the $m$-dimensional row vector ($e$ for 'entropy flow')
\bq
y_{e}=\frac{\partial K^c}{\partial p_S}
\eq
it follows that along the dynamics $X_K$ on $\cL$
\bq
\frac{d}{dt} S \geq y_{e}u
\eq
Hence $y_{e}$ is the output vector which is conjugate to $u$ in terms of {\it entropy flow}. The pair $(u, y_{e})$ is called the {\it flow of entropy} port of the system. 

All this is summarized in the following definition of a {\it port-thermodynamic system}.
\begin{df}[\cite{entropy}]
\label{def:portthermo}
Con-\\sider the manifold of extensive variables $\Q$ with coordinates $q^e=(E,S,q)$, and the cotangent bundle without zero section $\T^*\Q$ with coordinates $(q^e,p^e)=(E,S,q,p_E,p_S,p)$. A port-thermodynamic system on $\Q$ is a pair $(\cL,K)$, where $\cL \subset \T^*\Q$ is a Liouville submanifold describing the {\it state properties}, and $K=K^a + \sum_{j=1}^mK^c_ju_j,$ is a Hamiltonian on $\T^*\Q$, homogeneous of degree $1$ in $p^e$, and zero restricted to $\cL$, which generates the dynamics $X_K$.
Furthermore, $K^a$ is required to satisfy $\frac{\partial K^a}{\partial p_E}|_{\cL}=0$ and $\frac{\partial K^a}{\partial p_S}|_{\cL}\geq0$. The {\it power conjugate output} vector of the port-thermodynamic system is defined as $y_p=\frac{\partial K^c}{\partial p_E}$, and the {\it entropy flow conjugate output} vector as $y_{e}=\frac{\partial K^c}{\partial p_S}$.
\end{df}
Any port-thermodynamic system on $\T^*\Q$ projects to a corresponding contact system living on the projection of $\T^*\Q$ to the contact manifold $\mP(T^*\Q)$, where $\mP(T^*\Q)$ is the fiber bundle over $\Q$ with fibers the projective spaces $\mP(T_q^*\Q), q \in \Q$. In fact, see \cite{entropy,jgp} for details, since $\cL \subset \T^*\Q$ is a Liouville submanifold it projects to a Legendre submanifold $\widehat{\cL} \subset \mP(T^*\Q)$. Furthermore, since $K$ is homogeneous of degree $1$ in $p^e$ the homogeneous Hamiltonian vector field $X_K$ projects to a contact vector field $X_{\widehat{K}}$, with contact Hamiltonian $\widehat{K}$, that leaves $\widehat{\cL}$ invariant. Furthermore, by Euler's theorem both the power conjugate output $y_p$ and the entropy flow conjugate output $y_{e}$ are homogeneous of degree $0$, and thus project to functions on $\mP(T^*\Q)$. 

\medskip

Port-thermodynamic systems can be interconnected, either by their power ports or entropy flow ports, giving rise to a more complex port-thermodynamic system; cf. \cite{entropy} for details. For example, the power port interconnection of two port-thermodynamical systems with variables
\bq
(E_i,S_i,q_i,p_{E_i},p_{S_i},p_i) \in \T^*\Q_i, \quad i=1,2,
\eq
is defined as follows. With the homogeneity assumption in mind, impose the following constraint on the co-extensive variables
\bq
p_{E_1} = p_{E_2}=: p_E
\eq
This leads to the summation of the Liouville one-forms $\alpha_1$ and $\alpha_2$ given by
\bq
\alpha_{\mathrm{sum}}:= p_Ed(E_1 + E_2) + p_{S_1}dS_1 + p_{S_2}dS_2 + p_1dq_1 + p_2dq_2
\eq
on the {\it composed space} $\T^*\Q_1 \circ \T^*\Q_2$ defined as
\bq
\T^*\Q_1 \circ \T^*\Q_2 :=  \{(E,S_1,S_2,q_1,q_2, p_E, p_{S_1}, p_{S_2}, p_1,p_2)\} 
\eq
Let the state properties of the two systems be defined by the Liouville submanifolds $\cL_i \subset \T^*\Q_i, i=1,2$, then the state properties of the interconnected system are defined by the composition
\bq
\begin{array}{l}
\cL_1 \circ \cL_2 :=  \{(E,S_1,S_2,q_1,q_2, p_E, p_{S_1}, p_{S_2},p_1,p_2) \mid \\[2mm]
E=E_1 + E_2,  (E_i,S_i,q_i,p_{E_i},p_{S_i},p_i) \in \cL_i, \; i=1,2 \}
\end{array}
\eq
Furthermore, consider the dynamics on $\cL_i$ defined by Hamiltonians $K_i =K_i^a + K_i^cu_i, i=1,2$, where $K_i^c$ is the row vector of control Hamiltonians of system $i, i=1,2$. {\it Assume} that $K_i$ do {\it not} depend on the energy variables $E_i, i=1,2$. Then $K_1 +K_2$ is well-defined on $\cL_1 \circ \cL_2$ for all $u_1,u_2$. Next, consider the power conjugate outputs $y_{p1}, y_{p2}$. By imposing {\it interconnection constraints} on the power port variables $u_1,u_2,y_{p1},y_{p2}$ satisfying the {\it power preservation} property
\bq
y_{p1}^\top u_1 + y_{p2}^\top u_2=0,
\eq
then leads to an interconnected port-thermodynamic system with state properties described by $\cL_1 \circ \cL_2$. Similarly we can consider the entropy flow ports with entropy flow outputs $y_{e1},y_{e2}$, satisfying the {\it nonnegative entropy flow} property
\bq
y_{e1}^\top u_1 + y_{e2}^\top u_2 \geq 0,
\eq
leading again to a port-thermodynamic system.

\section{STABILITY ANALYSIS OF PORT-THERMODYNAMIC SYSTEMS}
\label{sec:asymp}

Let $(\cL,K)$ be a port-thermodynamic system, with $\cL \subset \T^*\Q$ a Liouville submanifold and $K=K^a + K^s$ a Hamiltonian which is homogeneous of degree $1$ in $p^e$ and zero on $\cL$, generating the dynamics $X_K$ on $\cL$. Consider the uncontrolled case $u=0$, so that $K^s$ (interaction with the surrounding) is zero. First we define the {\it equilibria} of port-thermodynamic systems. 
\begin{df}
The equilibria of $(\cL,K^a)$ are $(q^{e*},p^{e*}) \in \cL$ such that $X_{K^a}$ is zero at $(q^{e*},p^{e*})$, or equivalently $dK^a (q^{e*},p^{e*})=0$.
\end{df}
This is readily seen to be equivalent to the corresponding contact vector field $X_{\widehat{K}^a}$ being zero at the projection of $(q^{e*},p^{e*})$ to $\mP(T^*\Q)$, which is an element of the corresponding Legendre submanifold $\widehat{\cL}$.

How to assess the stability of equilibria of thermodynamic systems? Since the true dynamics is given by $X_{K^a}$ {\it restricted} to $\cL$ the stability of an equilibrium $(q^{e*},p^{e*}) = (E^*,S^*,q^*,p_E^*,p_S^*,p^*)$ is defined as the stability of $(q^{e*},p^{e*})$ with respect to the dynamics $X_{K^a}$ restricted to $\cL$. Consequently, contrary to standard Hamiltonian dynamics, the Hamiltonian $K^a$ is {\it not} a natural candidate Lyapunov function. 

\subsection{Stability analysis in the energy representation}
On the other hand, $\frac{\partial K^a}{\partial p_E}|_{\cL} =0$, or equivalently $\{K^a,E\}|_{\cL} =0$, where $E$ is the coordinate function on $\T^*\Q$ and $\{.,.\}$ is the standard Poisson bracket on $\T^*\Q$. Furthermore, in the energy representation of $\cL$
\bq
E=\bar{E}(S,q), \quad \mbox{ for all } (E,S,q,p_E,p_S,p) \in \cL,
\eq
and $\{K,E\}|_{\cL} =0$ is equivalent to 
\bq
\dot{\bar{E}} = \frac{\partial \bar{E}}{\partial S} \frac{\partial K^a}{\partial p_S} + \frac{\partial \bar{E}}{\partial q} \frac{\partial K^a}{\partial p} =0
\eq
Thus, whenever $\bar{E}$ has a {\it strict minimum} at $(S^*,q^*)$, then $\bar{E}$ is a Lyapunov function for the dynamics restricted to $\cL$, and stability of the equilibrium with respect to the dynamics on $\cL$ results by standard Lyapunov theory.

What can be done if $\bar{E}$ does not have a strict minimum at $(S^*,q^*)$ ? A classical tool in the stability analysis of ordinary Hamiltonian dynamics is to consider additional {\it conserved quantities}; see e.g. \cite{arnold, abraham,libermann}. In order to extend this idea to the present case let us strengthen our assumption on $K^a$ by requiring that $\frac{\partial K^a}{\partial p_E}=0$ everywhere on $\T^*\Q$; i.e., not just on $\cL$. Next we will consider additional conserved quantities for $X_{K^a}$ only depending on the extensive variables $S,q$; i.e., functions $C(S,q)$ such that
\bq
\{K^a,C\}=0
\eq
As a result also $\{K^a,E+ C \}=0$. Then we note the following identity regarding the Liouville form $\alpha$ on $\T^*\Q$
\bq
\begin{array}{l}
\alpha = p_EdE + p_SdS + pdq = \\[2mm]
 p_E d(E + C(S,q)) + (p_S - p_E\frac{\partial C}{\partial S}) dS + (p - p_E\frac{\partial C}{\partial q}) dq
\end{array}
\eq
Hence the transformation 
\bq
\begin{array}{l}
(E,S,q,p_E,p_S,p) \mapsto  \\
(E +C,S,q,p_E,p_S - p_E\frac{\partial C}{\partial S} , p - p_E\frac{\partial C}{\partial q}) \\
=:(\widetilde{E},\widetilde{S},\widetilde{q},\widetilde{p}_E,\widetilde{p}_S,\widetilde{p})
\end{array}
\eq
is a canonical point transformation (leaving the Liouville form $\alpha$ invariant). Note that in the new coordinates the intensive variables $\frac{p_S}{-p_E}, \frac{p}{-p_E}$ are transformed into {\it new} intensive variables 
\bq
\begin{array}{l}
\frac{\widetilde{p}_S}{-p_E} = \frac{p_S - p_E\frac{\partial C}{\partial S}}{-p_E} = \frac{p_S}{-p_E} + \frac{\partial C}{\partial S} \\[4mm]
\frac{\widetilde{p}}{-p_E} = \frac{p - p_E\frac{\partial C}{\partial q}}{-p_E} = \frac{p}{-p_E} + \frac{\partial C}{\partial q}
\end{array}
\eq
In these new coordinates the generating function for $\cL$ in entropy representation is given by $\bar{\widetilde{E}}(S,q)= \bar{E}(S,q) + C(S,q)$. Furthermore, since $\{K^a,E+ C\}=0$, the transformed Hamiltonian
\bq
\widetilde{K^a}(\widetilde{E},\widetilde{S},\widetilde{q},\widetilde{p}_E,\widetilde{p}_S,\widetilde{p}):=K^a(E,S,q,p_E,p_S,p)
\eq
satisfies $\{ \widetilde{K}, \widetilde{E} \}=0$. Hence in the new coordinates we are back to the situation considered before: if $\bar{\widetilde{E}}(S,q)$ has a strict minimum at $(S^*,q^*)$, then $\bar{\widetilde{E}}$ is a Lyapunov function for the dynamics restricted to $\cL$, and the equilibrium $(E^*,S^*,q^*)$ with $E^*=\bar{E}(S^*,q^*)$, is stable with respect to the dynamics on $\cL$.

Finally, note that the row vector $K^c$ of Hamiltonians in the new coordinates transforms to $\widetilde{K^c}(\widetilde{E},\widetilde{S},\widetilde{q},\widetilde{p}_E,\widetilde{p}_S,\widetilde{p})$, leading to the {\it transformed} power conjugate outputs
\bq
\widetilde{y}_p :=\frac{\partial \widetilde{K}^c}{\partial \widetilde{p}_E},
\eq
to be employed for asymptotic stabilization later on.

\subsection{Stability analysis in the entropy representation}
A similar analysis can be performed for the entropy representation of the Liouville submanifold $\cL$, where the entropy $S$ is expressed as a function $\bar{S}(E,q)$ of the energy variable $E$ and the other extensive variables $q$. By the Second Law we have $\frac{\partial K^a}{\partial p_S}|_{\cL} \geq 0$, or equivalently $\{K^a,S\}|_{\cL} \geq 0$. This could suggest to consider $-\bar{S}$ as a candidate Lyapunov function for the equilibrium $(E^*,q^*)$ at hand. However, typically the function $-\bar{S}(E,q)$ does {\it not} have a minimum at $(E^*,q^*)$. On the other hand, the entropy function $\bar{S}$ is known to be a {\it concave function} of $E,q$, and a classical idea is to consider instead the {\it availability function} $\bar{S}^*$ with respect to $(E^*,q^*)$, defined as
\bq
\begin{array}{l}
\A(E,q) := \bar{S}(E^*,q^*) -\bar{S}(E,q) \, + \\[2mm]
\frac{\partial \bar{S}}{\partial E}(E^*,q^*) (E - E^*) + \frac{\partial \bar{S}}{\partial q}(E^*,q^*) (q - q^*) 
\end{array}
\eq
(Note that $-\A$, but now as a function of $(E,q)$ and $(E^*,q^*)$, is also known as the Bregman divergence of $\bar{S}$.)
It is readily seen that the availability function  $\A$ is a {\it convex} function of $E,q$, which attains its {\it minimum} at $(E^*,q^*)$ with value zero. As a result, $\A^*$ is a candidate Lyapunov function for assessing the stability of $(E^*,q^*)$.

\begin{ex}
Chemical reaction networks are formulated as port-thermodynamic systems as follows; cf. \cite{louvain}. For simplicity we will not take volume and pressure into account, and consider concentrations instead of mole numbers. 
In the entropy representation $S$ is expressed as a function $\bar{S}=S(E,q)$, with $q \in \mR^m_+$ the vector of concentrations of the $m$ chemical species involved in the chemical reaction network. Recall that
\bq
\frac{\partial \bar{S}}{\partial q}(E,q)= - \frac{\mu}{T}, \quad \frac{\partial \bar{S}}{\partial E}(E,q)=\frac{1}{T}
\eq
The dynamics of the chemical reaction network is given by the homogeneous Hamiltonian vector field $X_{K^a}$, with
\bq
\begin{array}{rcl}
K^a &= & - p^\top Z \mL \Exp \frac{-Z^\top }{R}\frac{\partial S}{\partial q}(E,q) \, - \\[2mm]
&&p_S \frac{\partial S}{\partial q^\top }(E,q)Z \mL \Exp \frac{-Z^\top}{R}\frac{\partial S}{\partial q}(E,q)
\end{array}
\eq
Here $Z$ is the {\it complex composition matrix}, specifying the composition of the complexes (left- and right-hand sides of the chemical reactions) in the chemical species, $B$ is the {\it incidence matrix} of the graph of complexes, with edges corresponding to reactions, and $N=ZB$ is the stoichiometric matrix. $\Exp$ denotes the multi-dimensional exponential mapping; i.e., $\left(\Exp x\right)_i= \exp x_i$. Finally, $\mL = B^\top \mathcal{K} B$ is a weighted Laplacian matrix, with weights given by the diagonal elements of a certain matrix $\mathcal{K}$; see \cite{louvain} for details. 
Consider a thermodynamic equilibrium $(E^*,q^*)$, i.e., $B^\top Z^\top \mu^*=0$ where $\mu^*$ is the value of the chemical potential corresponding to $(E^*,q^*)$.
This implies that for an isolated chemical reaction network
\bq
\frac{d}{dt}S = \frac{1}{T} \mu^\top Z \mL \Exp (\frac{Z^\top \mu}{RT}) \geq 0,
\eq
with equality if and only if $B^TZ^\top \mu=N^\top \mu=0$, i.e., if and only if the {\it affinities} of the reactions are zero. Hence the equilibria of the system correspond to states of minimal (i.e., zero) entropy production. An equilibrium $(E^*,q^*)$ corresponds to By using $\A$ as Lyapunov function, it follows, under the standard assumption that trajectories will not converge to the boundary of the positive orthant $\mathbb{R}^m_+$, that any initial vector of concentrations in the positive orthant will converge to one of these equilibria; see e.g. the exposition in \cite{wang}. 
\end{ex}
In case other conserved quantities are needed we adopt the same strategy as in the energy representation, but now with $E$ and $S$ swapped. Thus we look for conserved quantities $C(E,q)$ such that $\{K^a,C\}=0$, and then consider the identity
\bq
\begin{array}{l}
\alpha = p_EdE + p_SdS + pdq = \\[2mm]
= (p_E - p_S\frac{\partial C}{\partial E}) dE + p_S d(S + C(E,q)) + (p - p_S\frac{\partial C}{\partial q}) dq
\end{array}
\eq
Similarly as before the mapping
\bq
\begin{array}{l}
(E,S,q,p_E,p_S,p) \mapsto  \\
(E,S +C,q,p_E - p_S\frac{\partial C}{\partial E}, p_S, p - p_S\frac{\partial C}{\partial q}) \\
=:(\widetilde{E},\widetilde{S},\widetilde{q},\widetilde{p}_E,\widetilde{p}_S,\widetilde{p})
\end{array}
\eq
defines a canonical point transformation, also leading to a modified entropy flow conjugate output $\widetilde{y}_{e}=\frac{\partial K^c}{\partial \widetilde{p}_S}$.

\subsection{Stability and asymptotic stabilization}\label{subsec:asymp}

Let $(E^*,S^*,q^*,p_E^*,p_S^*,p^*)$ be an equilibrium on $\cL$ as above. Consider the energy representation, and suppose that the function $\bar{E}(S,q)$, possibly with the help of an extra conserved quantity $C(S,q)$, has a strict minimum at $(S^*,q^*)$. By Lyapunov function theory it thus follows that the equilibrium is {\it stable}. How can we turn the equilibrium into an {\it asymptotically stable} equilibrium by feedback? The natural thing to do is to add {\it physical damping} by employing the power conjugate output $y_p$ (or, in case of an extra conserved quantity $C(S,q)$, the power conjugate output $\widetilde{y}_p$). How does this work for port-thermodynamic systems? 

Assume for simplicity of exposition that $m=1$ (scalar output $y_p$). Then consider an additional linear {\it damper system} (cf. \cite{entropy}), with Liouville submanifold
\bq
\cL_d= \{(U_d,S_d) \mid U_d=\bar{U}_d(S_d), p_{S_d}= -p_{U_d}\bar{U}_d'(S_d) \},
\eq
with entropy $S_d$ and internal energy $\bar{U}_d(S_d)$, with $\bar{U}_d'(S_d)$ its temperature. The dynamics of this damper system is generated by the Hamiltonian (see \cite{entropy},)
\bq
(p_{U_d} + p_{S_d} \frac{1}{U'(S_d)} ) du_d^2
\eq
(note the quadratic dependence on the input $u_d$), with power conjugate output $y_d=du_d$ (damping force). Then interconnect the system $(\cL,K=K^a + K^cu)$ to this damper system by setting
\bq
u= - y_d, \quad u_d= y_p
\eq
This results (after setting $p_{U_d}=p_E$) in the interconnected port-thermodynamic system with total Hamiltonian given as
\bq
\begin{array}{l}
K^a (E,S,q,p_E,p_S,p) - K^c(E,S,q,p_E,p_S,p)dy + \\[2mm]
(p_{U_d} + p_{S_d} \frac{1}{U'(S_d)} ) dy_p^2
\end{array}
\eq
with total energy $\bar{E}(S,q) + \bar{U}_d(S_d)$.
This implies that 
\bq
\frac{d}{dt} \bar{E}(S,q)= - \frac{d}{dt} \bar{U}_d(S_d) = - \bar{U}_d'(S_d) \dot{S_d}= - dy^2 \leq 0
\eq
Hence, by an application of LaSalle's Invariance principle, the system converges to the largest invariant set within the set where the power conjugate output $y_p$ is zero. Note that $y_p=0$ corresponds to zero entropy production $\dot{S_d}$; in accordance with irreversible thermodynamics (\cite{kondepudi}). If the largest invariant set where $y$ is zero equals the singleton $(E^*,S^*,q^*)$ then asymptotic stability of $(E^*,S^*,q^*)$ results; together with some limiting value $S_d^*$ of the entropy $S_d$ of the damper system. 

\section{STABILIZATION BY CONTROL BY INTERCONNECTION}\label{sec:int}
{\it Control by interconnection} is the paradigm of controlling a (plant) system by interconnecting it (through its inputs and outputs) to an additional controller system. The aim is to influence the dynamics of the original plant system by {\it shaping the dynamics} of the interconnected system by a proper choice of the controller system. 

Our treatment of asymptotic stabilization in the previous section is already an example of this methodology, since the given plant system was interconnected to a damper system, which can be considered to be the controller system. This methodology can be pursued much more generally; also invoking the use of conserved quantities as discussed before. As a very simple paradigmatic example, outside the normal thermodynamic realm, let us  consider the regulation of a mass-spring system to a non-zero set-point value of the spring extension and to zero velocity. This also serves a simple example of the theory of control by interconnection of port-Hamiltonian systems, cf. \cite{ortega, passivitybook, vds14}. We will show how to derive the same controller system within the more general framework of port-thermodynamic systems.

\begin{ex}
A mass-spring system with extensive variables $(E_p,z,\pi)$ is described by a Liouville submanifold $\cL$ of $\T^* \mR^3$, where $E_p$ is expressed as
\bq
E_p=\bar{E}_p(z,\pi) = \frac{1}{2}kz^2 + \frac{\pi^2}{2m},
\eq
with $z$ the extension of the spring with spring constant $k$, and $\pi$ the momentum of the mass with mass $m$. The dynamics is generated by the homogeneous Hamiltonian
\bq
K_p= \pi_z\frac{\pi}{m} - p_{\pi}kz + \left(p_{\pi} + p_E\frac{\pi}{m} \right) u_p,
\eq
where the input $u_p$ is the external force $u_p$ , and $y_p=\frac{\pi}{m}$ (velocity) is the power conjugate output. A scalar controller system with extensive variables $(E_c,q_c) \in \mR^2$ is given by the port-thermodynamic system $\cL_c,K_c)$, with energy $E_c$ expressed as $E_c=\bar{E}_c(q_c)$, and dynamics
\bq
K_c= \left(p_c + p_{E_c}\bar{E}_c'(q_c) \right)u_c
\eq
with output $y_c=\bar{E}_c'(q_c)$. The function $\bar{E}_c(q_c)$ is a design parameter, specifying the controller system.

The closed-loop system is obtained by the negative feedback (with $u$ the new input)
\bq
u_p=-y_c + u= -\bar{E}_c'(q_c) +u, \; u_c=y_p= \bar{E}_c'(q_c)
\eq
together with
\bq
E:=E_p+E_c,\; p_{E_p}=p_{E_c}=:p_E
\eq
This leads to the closed-loop Hamiltonian 
\bq
\begin{array}{rcl}
K & = & p_z\frac{\pi}{m} - p_{\pi}kz + \left(p_{\pi} + p_{E_p} \frac{\pi}{m}\right) \left(-\bar{E}'_c(q_c) +u \right) +\\[2mm]
&&+ \left(p_{q_c} + p_{E_c} \bar{E}'_c(q_c)\right) \frac{\pi}{m}
\end{array}
\eq
It is immediately seen that $C(z,\pi,q_c)= \Phi( z-q_c)$ for any function $\Phi: \mR \to \mR$ is a conserved quantity. This motivates to consider new canonical coordinates $(\widetilde{E}, \widetilde{z}, \widetilde{\pi}, \widetilde{q}_c, \widetilde{p}_E, \widetilde{p}_z, \widetilde{p}_{\pi}, \widetilde{p}_c)$, where 
\bq
\widetilde{E}=E+\Phi(w=z-q_c), \widetilde{p}_z= p_z - p_E \frac{\partial \Phi}{\partial w}, \widetilde{p}_{q_c} + p_E \frac{\partial \Phi}{\partial w},
\eq
while $\widetilde{z}=z, \widetilde{\pi}=\pi, \widetilde{q}_c=q_c,\widetilde{p}_E=p_E, \widetilde{p}_{\pi}=p_{\pi}$. In the new coordinates we compute $\widetilde{K}$ as
\bq
\begin{array}{l}
\widetilde{K}=\left( \widetilde{p}_z + \widetilde{p}_E \frac{\partial \Phi}{\partial w}\right)\frac{\widetilde{\pi}}{m} - \widetilde{p}_{\pi} k \widetilde{z} - \left(\widetilde{p}_{\pi} + \widetilde{p}_E\frac{\widetilde{\pi}}{m} \right) \bar{E}_c'(\widetilde{q}_c) + \\[2mm]
\left(\widetilde{p}_{q_c} - \widetilde{p}_E \frac{\partial \Phi}{\partial w} +
\widetilde{p}_E \bar{E}'_c(\widetilde{q}_c) \right) \frac{\widetilde{\pi}}{m} + 
\left(\widetilde{p}_{\pi} + \widetilde{p}_E \frac{\widetilde{\pi}}{m} \right) u ,
\end{array}
\eq
leading to the same power conjugate output $\widetilde{y}= y_p= \frac{\pi}{m}$ (velocity of the mass). \\
It is readily seen that for any $z^*$ the functions $\Phi$ and $\bar{E}_c$ can be chosen in such a way that the function $\bar{E} + \Phi$ has a strict minimum at the set-point value $(z^*,\pi^*=0, q_c^*)$ for some state value $q_c^*$ of the controller system. Furthermore, as shown in the Subsection \ref{subsec:asymp} this can be turned into asymptotic stabilization by further interconnecting the system with a damper system through the power port $(u,\widetilde{y}_p=y_p)$.
\end{ex}

\section{CONCLUSIONS}\label{sec:concl}
In this paper we have made initial steps towards a theory of 'control by interconnection' of port-thermodynamic systems. First we noted that candidate Lyapunov functions for an uncontrolled port-thermodynamic system can be inferred from the Liouville submanifold describing the state properties: either the energy expressed as a function of the other extensive variables, or (minus) the entropy as a function of the remaining extensive variables. Furthermore, we can employ conserved quantities for the dynamics in order to shape these functions, as reflected in new canonical coordinates for the total space of extensive and co-extensive variables. 
This should be contrasted with control by interconnection of port-Hamiltonian systems, where a Lyapunov function is sought to be constructed on the basis of the Hamiltonian, as well as conserved quantities. As a first example of control by interconnection it was shown how a stable equilibrium of a port-thermodynamic system can be rendered asymptotically stable by the interconnection with a damper system, having its own energy and entropy. Finally, in Section \ref{sec:int} it was shown how set-point regulation for a simple mass-spring system can be achieved by the use of two additional port-thermodynamic systems; one for energy shaping and one for asymptotic stabilization. Future research is concerned with the extension of this methodology to non-trivial thermodynamic situations, such as the Continuous Stirred Tank Reactor.

\end{document}